\documentclass[11pt, fleqn]{amsart}
\usepackage{amssymb,amsmath,latexsym,amsbsy,color,enumerate}
\numberwithin{equation}{section}
\def\NN{\mbox{$I\hspace{-.06in}N$}}

\def\RR{\mbox{$I\hspace{-.06in}R$}}
\def\CC{\mbox{$C\hspace{-.11in}\protect\raisebox{.5ex}{\tiny$/$}
\hspace{.06in}$}}

\newtheorem{theorem}{Theorem}

\newtheorem{lemma}{Lemma}
\newtheorem{remark}{Remark}

\begin{document}
 \title{Error of Tikhonov's regularization for  integral convolution equations}
    \author{Dang Duc Trong}
    \address{Dang Duc Trong, Department of Mathematics, Hochiminh City National University, 227 Nguyen Van Cu, Q5, HoChiMinh City, Vietnam}
    \email{ddtrong@mathdep.hcmuns.edu.vn}
\author{Truong Trung Tuyen}
    \address{Truong Trung Tuyen, Department of Mathematics, Indiana University Bloomington,  Rawles Hall, IN 47405}
 \email{truongt@indiana.edu}
    \date{\today}
    \begin{abstract}
Let $\varphi$ be a nontrivial function of $L^1(\RR )$. For each $s\geq 0$ we put
\begin{eqnarray*}
p(s)=-\log \int _{|t|\geq s}|\varphi (t)|dt.
\end{eqnarray*}
If $\varphi$ satisfies
\begin{equation}
\lim _{s\rightarrow \infty}\frac{p(s)}{s}=\infty ,\label{170506.1}
\end{equation}
we obtain asymptotic estimates of the size of small-valued sets $B_{\epsilon}=\{x\in\RR :~|\widehat{\varphi}(x)|\leq \epsilon ,~|x|\leq R_{\epsilon}\}$ of Fourier transform
\begin{eqnarray*}
\widehat{\varphi}(x)=\int _{-\infty}^{\infty}e^{-ixt}\varphi (t)dt,~x\in \RR ,
\end{eqnarray*}
in terms of $p(s)$ or in terms of its Young dual function
\begin{eqnarray*}
p^{*}(t)=\sup _{s\geq 0}[st-p(s)],~t\geq 0.
\end{eqnarray*}
Applying these results, we give an explicit estimate for the error of Tikhonov's regularization for the solution $f$ of the integral convolution equation
\begin{eqnarray*}
\int _{-\infty}^{\infty}f(t-s )\varphi (s )ds =g(t),
\end{eqnarray*}
where $f,g \in L^2(\RR )$ and $\varphi$ is a nontrivial function of $L^1(\RR )$ satisfying condition (\ref{170506.1}), and $g,\varphi$ are known non-exactly. Also, our results extend some results of \cite{tld} and \cite{tqd}.
\end{abstract}
    \maketitle
\section{Introduction}
Let $\varphi$ be a nontrivial function of $L^1(\RR )$. Then its Fourier transform is defined by
\begin{eqnarray*}
\widehat{\varphi}(x)=\int _{-\infty}^{\infty}e^{-ixt}\varphi (t)dt,~x\in \RR .
\end{eqnarray*}
If $\varphi$ is of compact support, there associates an entire function of exponential type of Laplace transform type
\begin{equation}
\Phi (z)=\int _{-\infty}^{\infty}e^{zt}\varphi (t)dt.\label{fou}
\end{equation}
Generally, if $\varphi\in L^1(\RR )$  satisfies
\begin{equation}
\int _{-\infty}^{\infty}e^{|t|s}|\varphi (t)|dt <\infty , \label{150506}
\end{equation}
for all $s>0$ then $\Phi$ in (\ref{fou}) is defined on all over the complex plane $\CC$ and is an entire function, but its order may be any positive number. If $\varphi \geq 0$ then condition (\ref{150506}) is also necessary for $\Phi$ to be an entire function. Indeed, in this case we have, for every $R>0$,
\begin{eqnarray*}
\max _{|z|\leq R}|\Phi (z)|\geq \frac{1}{2}(\Phi (R)+\Phi (-R))\geq \frac{1}{2}\int _{-\infty}^{\infty}e^{|t|R}\varphi (t)dt.
\end{eqnarray*}

As we will show in Section 2, the growth of $\Phi (z)$ is related to the function
\begin{equation}
p(s)=-\log \int _{|t|\geq s}|\varphi (t)|dt ,~s\geq 0.\label{170506.2}
\end{equation}
In fact, Theorem \ref{theo3} show that Young dual function $p^{*}(s)$ of $p(t)$, defined by
\begin{eqnarray*}
p^{*}(s)=\sup _{t\geq 0}[st-p(t)],~s\geq 0,
\end{eqnarray*}
is an appropriate quantity to estimate the growth of $\log |\Phi (z)|$. It turns out that $\varphi$ satisfies (\ref{150506}) if and only if it satisfies
\begin{equation}
\lim _{s\rightarrow\infty}\frac{p(s)}{s}=\infty .\label{150506.1}
\end{equation}

Properties of Fourier and Laplace transforms are of great interest because these functions often occur in different areas of calculus. As a typical application, Fourier transforms give us the direct solution of integral convolution equations. These problems, so-called deconvolution problems (see, e.g.\cite{ram}), arise in applications: an input signal $f$ generates an output signal $g$ by the formula
\begin{eqnarray*}
Af=\int _D\varphi (x,\gamma )f(\gamma )d\gamma =g(x),~x\in D.
\end{eqnarray*}
Often one has
\begin{equation}
Af=\int _{-\infty}^{\infty}\varphi (t-s)f(s)ds=g(t),~t\in\RR ,\label{i1}
\end{equation}
where $\varphi (t)\in L^1(\RR )$. The deconvolution problem consists of finding $f$, given $g$ and $\varphi (x,\gamma )$. For (\ref{i1}),  another problem, so-called the identification problem, is of practical interest: given $f$ and $g$, find $\varphi (t)$. The function $\varphi (t)$ characterizes the linear system which generates the output $g(t)$ given the input $f(t)$.

Now if $f$ is the solution of (\ref{i1}), applying Fourier transform to both sides of (\ref{i1}) we get
\begin{equation}
\widehat{f}=\frac{\widehat{g}}{\widehat{\varphi}}.\label{eq11}
\end{equation}
Then apply the inverse Fourier transform we find $f$. We may see easily from equation (\ref{eq11}) that $f$ depends nonlinearly on $\varphi$.

As easily seen, in case $\varphi $ has compact support, the set of zeros of $\widehat{\varphi}(z)\equiv \Phi (-iz)$ effects largely the recovering the function $f$ from its Fourier transform. Since $\widehat{\varphi}(z)$ is a function of class Cartwright, that is
\begin{eqnarray*}
\int _{-\infty}^{\infty}\frac{\log ^+|\widehat{\varphi}(t)|}{1+t^2}dt\leq \int _{-\infty}^{\infty}\frac{\log ^+||\varphi||_{L^1(\RR )}}{1+t^2}dt<\infty ,
\end{eqnarray*}
the distribution of its zeros is well known (see, e.g., \cite{lev}). In a recent paper \cite{sed}, Sedletskii obtained some detailed properties of the set of zeros of $\Phi$ for a subclass of functions $f$ that is of $L^1(0,1)$ and positive (see also references in \cite{sed}).

However, for really computing the solution $f$ and for the regularization of equation (\ref{i1}), we must know more about the structure of the small-valued sets of $\Phi$. This problem goes back to the well-known theorem of Cartan about the size of the small-valued sets $B_{\epsilon}=\{z\in \CC :~|P(z)|\leq \epsilon\}$ where $P(z)$ is a polynomial. He proved that $B_{\epsilon}$ is contained in a finite of disks whose sum of radii is less than $C_1\epsilon ^{\frac{1}{n}}$ where $n$ is the degree of $P(z)$ and $C_1$ is a constant that depends only on the leading coefficient of $P(z)$ and $n$ (see Theorem 3 of $\S 11.2$ in \cite{lev}). In particularly we have
\begin{equation}
\lim _{\epsilon\rightarrow 0}m(B_{\epsilon})=0, \label{i11}
\end{equation}
where $m(.)$ is the Lebesgue's measure. For general case of entire functions, the conclusion (\ref{i11}) does not hold, for even simple functions such as $e^z$, instead of many results representing the structure of $B_{\epsilon}$ (see, e.g., \cite{lev}).

Recently, in seeking regularization schemes for the problem of determination of heat source in one and two dimensional, the authors in \cite{tld} and \cite{tqd} proved that $\Phi $ in (\ref{fou}) has only finitely many zeros on the positive real-axis if  $\varphi$ is of a fairly wide subclass of $L^2(0,1)$. Indeed, as easily seen, the results in \cite{tld} and \cite{tqd} depend only on the estimate of the size of the set
\begin{eqnarray*}
B_{\epsilon}=\{x\in \RR :~|\Phi (x)|\leq \epsilon , |x|\leq R_{\epsilon}\},
\end{eqnarray*}
where $\gamma ,\Phi$ are as above and $R_{\gamma}>0$ depends on $\gamma$.

In this paper, we will give some asymptotic estimates for the small-valued sets
\begin{eqnarray*}
B_{\epsilon}=\{x\in \RR :~|\widehat{\varphi} (x)|\leq \epsilon , |x|\leq R_{\epsilon}\},
\end{eqnarray*}
to functions $\varphi$ satisfying condition (\ref{150506.1}) and apply these estimates to the Tikhonov's regularization of Problem (\ref{i1}). Also, our results extend some results in \cite{tld} and \cite{tqd}. For convenience, we recall briefly Tikhonov's regularization before stating results.

The difficulty of applying direct formula (\ref{eq11}) arises in two aspects: one, in reality we can not have the exact data $(\varphi _0,g_0)$ but only the measured data $(\varphi _{\epsilon},g_{\epsilon })$, and two, even for the case in which we have the exact data, the inverse Fourier transform may not be efficiently computed if $\widehat{\varphi _0}$ has zeros on the real axis. Hence a regularization is needed.

Let $\epsilon >0$, we denote by $(\varphi _0,g_0)$ the exact data and $(\varphi _{\epsilon},g_{\epsilon})$ the measured data with the tolerable error  $\epsilon$, that is
\begin{equation}
\epsilon \geq\max \{||\varphi _0-\varphi _{\epsilon}||_{L^1(\RR)},||g_0-g_{\epsilon}||_{L^2(\RR)}\}.\label{170506.3}
\end{equation}

Tikhonov's regularization applied to equation (\ref{i1}) is to construct an approximation $f_{\epsilon}$ whose Fourier transform is
\begin{equation}
\widehat{f_{\epsilon}}=\frac{\widehat{g_{\epsilon}}\overline{\widehat{\varphi _{\epsilon}}}}{\delta _{\epsilon}+|\widehat{\varphi _{\epsilon}}|^2},\label{170506.4}
\end{equation}
where $\delta _{\epsilon}>0$ is a regularization parameter appropriately chosen depending on $\epsilon$.

We have the following general estimate for the error $||f_0-f_{\epsilon}||_{L^2(\RR )}$ (see also \cite{are} and Chapter 10 in \cite{bau} for similar estimates).

\begin{theorem}
Let $\beta$ be a constant in $(0,\frac{1}{3})$. Let $f_0\in L^2(\RR )$ be the exact solution of (\ref{i1}) corresponding to $\varphi _0\in L^1(\RR ),~g_0\in L^2(\RR )$ . Let $\epsilon >0$ be small, $\varphi _{\epsilon}, g_{\epsilon}$ be as in (\ref
{170506.3}) and let $f_{\epsilon}$ be as in (\ref
{170506.4}), where $\delta _{\epsilon}=\sqrt[4]{\frac{C_1}{C_2}}\epsilon ^{(1+3\beta )/2}$ and
\begin{eqnarray*}
C_1&=&4(1+||g_0||^2_{L^2(\RR )}+||\varphi _0||^2_{L^1(\RR )}),\\
C_2&=&1+||g_0||^2_{L^2(\RR )}.
\end{eqnarray*}

Let $R_{\epsilon }>0$ be a sequence such that
\begin{eqnarray*}
\lim _{\epsilon\rightarrow 0}R_{\epsilon}=\infty .
\end{eqnarray*}
 Then
\begin{eqnarray*}
||f_0-f_{\epsilon}||^2_{L^2(\RR )}&\leq& 3[\int _{|\lambda |>R_{\epsilon},|\widehat{\varphi _0} (\lambda )|<\epsilon ^{\beta}}|F(f_0)|^2d\lambda \\
&+&\int _{|\lambda |<R_{\epsilon},|\widehat{\varphi _0} (\lambda )|<\epsilon ^{\beta}}|F(f_0)|^2d\lambda + 2\sqrt{C_1C_2}\epsilon ^{1-3\beta}],
\end{eqnarray*}
where $F(f_0)(\lambda )=\int _{-\infty}^{\infty}f_0(x)e^{-i\lambda x}dx$ is the Fourier transform of $f_0$ and $\widehat{\varphi _0}(\lambda )=\int _{-\infty}^{\infty}\varphi _0(x)e^{-i\lambda x}dx$ is the Fourier transform of $\varphi _0$.
\label{Tikhonov}\end{theorem}

If $\varphi _0$ satisfies condition (\ref{150506}), we can show that the bound of $||f_0-f_{\epsilon}||_{L^2(\RR )}$ in above Theorem decreases to zero as $\epsilon$ decreases to zero. This assertion comes easily from
\begin{eqnarray*}
&&\lim _{\epsilon\rightarrow 0}\int _{|\lambda |>R_{\epsilon},|\widehat{\varphi_0} (\lambda )|<\epsilon ^{\beta}}|F(f_0)|^2d\lambda +\int _{|\lambda |<R_{\epsilon},|\widehat{\varphi_0} (\lambda )|<\epsilon ^{\beta}}|F(f_0)|^2d\lambda\\
&=&\lim _{\epsilon\rightarrow 0}\int _{\RR}\chi _{\{|\widehat{\varphi_0} (\lambda )|<\epsilon ^{\beta}\}}|F(f_0)|^2d\lambda =0,
\end{eqnarray*}
by Lebesgue's dominated convergence theorem.

However, to get a more explicit estimate for $||f_0-f_{\epsilon}||_{L^2(\RR )}$ some a-priopri information about $f_0$ must be assumed. For example, in Chapter 10 of \cite{bau} (see also \cite{are}), the author gave estimates for $||f_0-f_{\epsilon}||_{L^2(\RR )}$ in the case the kernel $\varphi _0$ is of two special types of $L^1$ functions (see conditions (10.4) and (10.10) in \cite{bau}) and in additionally, the solution $f_0$ satisfies the condition
\begin{equation}
|\widehat{f_0}(\lambda )|^2\leq C_1(1+|\lambda |^2)^{-q},\lambda \in \RR ,\label{240406}
\end{equation}
for some $q>\frac{1}{2}$.

Condition (10.4) in \cite{bau} is a set of four conditions and condition (10.10) in \cite{bau} is that
\begin{equation}
|\widehat{\varphi _0}(\lambda )|^2\geq C_0e^{-a\lambda },~\lambda >0,\label{240406.1}
\end{equation}
where $a,C>0$ are constants. While condition (\ref{240406}) imposed on $f_0$ is natural and is easily verified (for example, if $f_0\in W^{1,1}(\RR )$ we have $|\widehat{f_0}(\lambda )|^2\leq C_0(1+|\lambda |^2)^{-1}$), conditions (10.4) and (10.10) above may not hold for a general kernel $\varphi _0\in L^1(\RR )$ and fairly difficult to verify for a concrete kernel $\varphi _0$ in particularly in the case $\varphi _0$ is  known non-exactly.

In this paper, we will give estimates for error of Tikhonov's regularization for Problem (\ref{i1}) for the case the solution $f_0$ satisfies (\ref{240406}) and the kernel $\varphi _0$ satisfies (\ref{150506.1}). As we will show in Section 2, this class consists of functions satisfying condition (\ref{150506}).

For each $\epsilon >0$ we put
\begin{equation}
s_{\epsilon}=\inf \{s>0 :~e^{-p(s)} \leq \epsilon\},\label{150506.2}
\end{equation}
where $p(s)$ is as in (\ref{150506.1}).

Our main result is the following theorem. The estimate of our result in the case $\varphi _0$ has compact support is comparable to the estimate $||f_0-f_{\epsilon }||_{L^2(\RR )}\leq C(\log \frac{1}{\epsilon} )^{-q+\frac{1}{2}}$ gave in Theorem 10.8 in \cite{bau} (which uses additionally condition (\ref{240406.1}) in its proof).

\begin{theorem}
Let assumptions be as in Theorem \ref{Tikhonov}. In additional, assume that $f_0$ satisfies condition (\ref{240406}) and $\varphi _0$ satisfies condition (\ref{150506.1}). For $\epsilon >0$ small enough choose $s_{\epsilon}$ as in (\ref{150506.2}) and $R_{\epsilon} >0$ satisfying
\begin{eqnarray*}
[(q+\frac{1}{2})\log R_{\epsilon}+\log (15e^3)][\log ||\varphi _0||_{L^1(\RR )}+2es_{\epsilon}R_{\epsilon}]=-\log (\epsilon ^{\beta}+\epsilon).
\end{eqnarray*}
Then
\begin{eqnarray*}
\lim _{\epsilon\rightarrow 0}R_{\epsilon}=\infty ,
\end{eqnarray*}
and there exists a constant $C_3>0$ independent of $\epsilon$ such that
\begin{eqnarray*}
||f_0-f_{\epsilon}||_{L^2(\RR )}\leq C_3R_{\epsilon}^{-q+\frac{1}{2}},
\end{eqnarray*}
for $\epsilon >0$ small enough. In particular, in case $\varphi _0$ is of compact support then $s_{\epsilon}$ is bounded for $\epsilon >0$ small enough, hence
\begin{eqnarray*}
\lim _{\epsilon\rightarrow 0}\frac{\log R_{\epsilon}}{\log\log\frac{1}{\epsilon}}=1.
\end{eqnarray*}
\label{fgood}\end{theorem}
\begin{proof}
The existence of $R_{\epsilon}$ is easily seen by considering the function
\begin{eqnarray*}
f(R)=[(q+\frac{1}{2})\log R+\log (15e^3)][\log ||\varphi _0||_{L^1(\RR )}+2es_{\epsilon}R],~R\geq 0.
\end{eqnarray*}
Condition (\ref{150506.1}) gives
\begin{eqnarray*}
\lim_{\epsilon\rightarrow 0}R_{\epsilon}=\infty .
\end{eqnarray*}
By Theorem \ref{Tikhonov} and condition (\ref{240406}), there exists $C_4>0$ such that
\begin{eqnarray*}
||f_0-f_{\epsilon}||_{L^2(\RR )}\leq C_4(R^{-q+\frac{1}{2}}+m(B_{\epsilon})+\epsilon ^{1-3\beta}),
\end{eqnarray*}
where $m(B_{\epsilon})$ is the Lebesgue's measure of the set
\begin{eqnarray*}
B_{\epsilon}=\{x\in\RR :~|\widehat{\varphi _0}(x)|\leq \epsilon ^{\beta},~|x|\leq R_{\epsilon}\}.
\end{eqnarray*}
Now apply Theorem \ref{theo2} we get the conclusion of Theorem \ref{fgood}.
\end{proof}

The rest of this paper consists of three sections. In Section 2 we explore some properties of entire functions of Laplace transform type. In Section 3 we
state and prove estimates of the size of small-valued sets for entire functions of Laplace transform type and for Fourier transforms. In Section 4 we
prove Theorem \ref{Tikhonov}.
\section{Properties of Laplace transforms}
In this section we explore some properties of entire functions of type (\ref{fou}), where $\varphi\in L^1(\RR )$ satisfies condition (\ref{150506}).

First we consider the case $\varphi$ has a compact support. Without loss of  generality, we assume the support of $\varphi$ is contained in $[0,1]$. So the function $\Phi$ in (\ref{fou}) take a simpler form
\begin{eqnarray*}
\Phi (z)=\int _{0}^{1}e^{zt}\varphi (t)dt,
\end{eqnarray*}
where $\varphi (t)$ is a nontrivial function of $L^1(0,1)$.

We put
\begin{eqnarray*}
\sigma =\inf \{a\in [0,1]: \varphi |_{[a,1]}=0~a.e.\},\\
\mu =\sup \{a\in [0,1]:\varphi |_{[0,a]}=0~a.e.\}.
\end{eqnarray*}

\begin{lemma} (On order and type of $\Phi$)

i) Order of $\Phi$ is unity.

ii) Type of $\Phi$ is $\sigma$.

iii) We have
\begin{eqnarray*}
\limsup _{R\rightarrow \infty}\frac{\log |\Phi (R)|}{R}&=&\sigma ,\\
\limsup _{R\rightarrow \infty}\frac{\log |\Phi (-R)|}{R}&=&-\mu .
\end{eqnarray*}
\label{lem1}\end{lemma}
Before proving Lemma \ref{lem1} we recall that (see e.g. \cite{lev}) if $\Phi (z)$ is an entire function then its order $\rho$ (or $ord (\Phi )$ for short) and type $\sigma$ are defined by
\begin{eqnarray*}
\rho &=&\limsup_{r\rightarrow\infty}\frac{\log\log M_{\Phi}(r)}{\log r},\\
\sigma&=&\limsup _{r\rightarrow\infty}\frac{\log M_{\Phi}(r)}{r^{\rho}},
\end{eqnarray*}
where
\begin{eqnarray*}
M_{\Phi (r)}=\max _{|z|\leq r}|\Phi (z)|.
\end{eqnarray*}
\begin{proof} i) We have
\begin{eqnarray*}
|\Phi (z)|\leq Ce^{|z|},
\end{eqnarray*}
for $z\in\CC$ and $C=||\varphi||_{L^1(\RR )}$. It follows that $ord(\Phi )\leq 1$.

Now we show $ord (\Phi )=1$ . Assume by contradiction that $ord (\Phi )<1$.

Applying Theorem 1 in $\S 6.1$ in \cite{lev} for domains $\{-\frac{\pi}{2}<\arg z<\frac{\pi}{2}\}$ and $\{\frac{\pi}{2}< \arg z<\frac{3\pi}{2}\}$, noting that $\Phi$ is bounded on the imaginery-axis, we have that $\Phi$ is bounded on $\CC$. Thus $\Phi$ is a constant, so $\Phi '(z)\equiv 0$. Thus we have $t\varphi\equiv 0$ or equivalently $\varphi \equiv 0$, since
\begin{eqnarray*}
\Phi '(z)=\int _{0}^{1}e^{zt}t\varphi (t)dt,
\end{eqnarray*}
which is a contradiction.

ii) We have the type of $\Phi$ is less than or equal to $\sigma$ because
\begin{eqnarray*}
|\Phi (z)|=|\int _0^{\sigma}e^{zt}\varphi (t)dt|\leq e^{\sigma |z|}\int _{0}^{\sigma}|\varphi (t)|dt.
\end{eqnarray*}

We prove an equivalent result: If $a\in [0,1]$ is so that $\varphi$ is not identity to zero in $[a,1]$ then the type of $\Phi$ is greater than or equal to $a$.

We assume by contradiction that the type of $\Phi$ is less than $a$.

We consider the function
\begin{eqnarray*}
\Gamma (z)=e^{-az}\int _{a}^{1}e^{zt}\varphi (t)dt.
\end{eqnarray*}

If $z>0$ we have
\begin{eqnarray*}
|\Gamma (z)|&=&|e^{-az}\Phi (z)-e^{-az}\int _{0}^{a}e^{zt}\varphi (t)|\\
&\rightarrow& 0
\end{eqnarray*}
when $z\rightarrow +\infty$ on $\RR$. Indeed, because of our contradiction assumption that the type of $\Phi$ is less than $a$ it follows
\begin{eqnarray*}
\lim _{z\rightarrow\infty}e^{-az}\Phi (z)=0,
\end{eqnarray*}
while by Lebesgue's dominated convergence theorem applied to the sequence $g_z(t)=e^{z(t-a)}\varphi (t)$ on $[0,a]$ it follows that
\begin{eqnarray*}
\lim _{z\rightarrow\infty}e^{-az}\int _0^{a}e^{zt}\varphi (t)dt=\lim _{z\rightarrow\infty}\int _0^{a}g_z(t)dt=0.
\end{eqnarray*}

If $z<0$, by Lebesgue's dominated convergence theorem as above we get
\begin{eqnarray*}
|\Gamma (z)|&=&e^{a|z|}|\int _{a}^{1}e^{-|z|t}\varphi (t)dt|\\
&\rightarrow& 0
\end{eqnarray*}
when $z\rightarrow -\infty$ on $\RR$.

If $z$ is purely imaginery, we have
\begin{eqnarray*}
|\Gamma (z)|\leq C.
\end{eqnarray*}

Applying Theorem 1 in $\S 6.1$ in \cite{lev} for domains $\{0<\arg z<\frac{\pi}{2}\}$, $\{\frac{\pi}{2}<\arg z<\pi \}$, $\{\pi <\arg z<\frac{3\pi}{2}\}$ and $\{\frac{3\pi}{2}<\arg z<2\pi \}$ we get $|\Gamma (z)|\leq C$ for all $z\in\CC$, which implies that $\Gamma (z)=$constant, a contradiction (qed).

iii) We prove similarily to the proof of ii). We prove first that
\begin{eqnarray*}
\limsup _{R\rightarrow \infty}\frac{\log |\Phi (-R)|}{R}=-\mu .
\end{eqnarray*}

If $R>0$ we have
\begin{eqnarray*}
|\Phi (-R)|=|\int _{\mu }^{1}e^{-Rt}\varphi (t)dt|\leq e^{-R\mu}||\varphi ||_{L^1(0,1)},
\end{eqnarray*}
hence
\begin{eqnarray*}
\limsup _{R\rightarrow\infty}\frac{\log |\Phi (-R)|}{R}\leq -\mu .
\end{eqnarray*}

Thus, in order to prove the equality required, we need to show only that: if $a>\mu $ then
\begin{eqnarray*}
\limsup _{R\rightarrow\infty}\frac{\log |\Phi (-R)|}{R}\geq -a.
\end{eqnarray*}

Assume by contradiction that
\begin{eqnarray*}
\limsup _{R\rightarrow\infty}\frac{\log |\Phi (-R)|}{R}< -a,
\end{eqnarray*}
for some $a>\mu$. We consider
\begin{eqnarray*}
\Gamma (z)=e^{-az}\int _0^{a}e^{zt}\varphi (t)dt.
\end{eqnarray*}

By arguments in ii) we have
\begin{eqnarray*}
\lim _{R\rightarrow \infty}|\Gamma (R)|&=&0,\\
\lim _{R\rightarrow \infty}|\Gamma (-R)|&=&\lim _{R\rightarrow \infty}|e^{aR}\Phi (-R)-e^{aR}\int _a^{1}e^{-Rt}\varphi (t)dt|=0.
\end{eqnarray*}

Thus as in ii) we have $\Gamma (z)$ is a constant, which is a contradiction.

The proof of
\begin{eqnarray*}
\limsup _{R\rightarrow \infty}\frac{\log |\Phi (R)|}{R}=\sigma
\end{eqnarray*}
is similar and indeed was contained in the proof of ii).
\end{proof}

\begin{lemma} (On zeros of $\Phi$)
\begin{eqnarray*}
\lim _{R\rightarrow\infty}n(R)=\infty
\end{eqnarray*}
and
\begin{eqnarray*}
\lim _{R\rightarrow\infty}\frac{n(R)}{R}=\frac{d}{\pi}:= \frac{\sigma -\mu}{\pi} ,
\end{eqnarray*}
where $n(R)=\{z:~|z|\leq R,~\Phi (z)=0\}$.
\label{lem3}\end{lemma}

\begin{proof} First we show that $\Phi$ has infinitely many zeros. Indeed, if $\Phi$ has only finitely many zeros, since $\Phi$ is of order $1$, by Hadamard's theorem we may write
\begin{eqnarray*}
\Phi (z)=e^{az}P(z),
\end{eqnarray*}
where $P$ is a polynomial.

So if we differentiate $e^{-az}\Phi (z)$ $m$ times, where $m$ is greater than the order of $P(z)$ we obtain
\begin{eqnarray*}
\int _{0}^{1}e^{z(t-a)}(t-a)^m\varphi (t)dt\equiv 0,
\end{eqnarray*}
so $\varphi \equiv 0$, a contradiction.

We have $|\Phi (it)|\leq C=||\varphi ||_{L^1(\RR )}$ if $t\in \RR$. So
\begin{eqnarray*}
\int _{-\infty}^{\infty}\frac{\log ^+|\Phi (it)|}{1+t^2}dt \leq \int _{-\infty}^{\infty}\frac{\log ^+C}{1+t^2}dt<\infty .
\end{eqnarray*}
Thus $\Phi (iz)$ is of class Cartwright (see Lecture 16 in \cite{lev}).

Applying Theorem 1 of $\S 17.2$ in \cite{lev}, we have
\begin{eqnarray*}
\lim _{R\rightarrow\infty}\frac{n(R)}{R}=\frac{d}{\pi}
\end{eqnarray*}
where $d$ is the width of the indicator diagram of $\Phi (iz)$, i.e.
\begin{eqnarray*}
d=\limsup _{R\rightarrow\infty}\frac{\log |\Phi (R)|}{R}+\limsup _{R\rightarrow\infty}\frac{\log |\Phi (-R)|}{R}.
\end{eqnarray*}

By Lemma \ref{lem1} iii) we have
\begin{eqnarray*}
\limsup _{R\rightarrow\infty}\frac{\log |\Phi (R)|}{R}&=&\sigma ,\\
\limsup _{R\rightarrow\infty}\frac{\log |\Phi (-R)|}{R}&=& -\mu .
\end{eqnarray*}

So $d= \sigma -\mu$.
\end{proof}

\begin{lemma} (On representation of $\Phi$).
Denote $\CC ^+=\{z\in \CC:~Im z>0\}$. i.e. $\CC ^+$ is the upper half-plane. Then
\begin{eqnarray*}
\Phi (z)=Cz^me^{\frac{\sigma +\mu}{2} z}\prod _{z_i\in \RR}(1-\frac{z}{z_i})\prod _{z_i\in \CC^+}(1-z(\frac{1}{z_i}+\frac{1}{\overline{z_i}})+\frac{z^2}{|z_i|^2}),
\end{eqnarray*}
where $C\in \RR$.
\label{lem2}\end{lemma}

\begin{proof}  Applying Theorem 1 of $\S 17.2$ in \cite{lev} for $\Phi (iz)$, if $\{z_k\},z_k\not=0$ are zeros of $\Phi$ then
\begin{eqnarray*}
\sum _{k}|Re \frac{1}{z_k}|<\infty .
\end{eqnarray*}

Applying Hadamard's theorem, noting that $\Phi (\overline{z})=\overline{\Phi (z)}$, we may write
\begin{eqnarray*}
\Phi (z)=Cz^me^{az}\lim _{r\rightarrow\infty}[\prod _{z_i\in \RR ,|z_i|\leq r}(1-\frac{z}{z_i})\prod _{z_i\in \CC ^+ ,|z_i|\leq r}(1-\frac{z}{z_i})(1-\frac{z}{\overline{z_i}})e^{z(\frac{1}{z_i}+\frac{1}{\overline{z_i}})}].
\end{eqnarray*}

Since $ord (\Phi) =1$ we have
\begin{eqnarray*}
\sum _{z_i}\frac{1}{|z_i|^2}<\infty .
\end{eqnarray*}

This, combined with the event that $\frac{1}{z_i}+\frac{1}{\overline{z_i}}=2Re \frac{1}{z_i}$ and the inequality above
\begin{eqnarray*}
\sum _{k}|Re \frac{1}{z_k}|<\infty ,
\end{eqnarray*}
allows us to write
\begin{eqnarray*}
\Phi (z)=Cz^me^{bz}\prod _{z_i\in \RR}(1-\frac{z}{z_i})\prod _{z_i\in \CC ^+ }(1-\frac{z}{z_i})(1-\frac{z}{\overline{z_i}}).
\end{eqnarray*}

That $C,b\in \RR$  is easy to see. Now we show that $b= \frac{\sigma +\mu}{2}$.

If $z_i\notin \RR$  and $z\in\CC$ we have
\begin{eqnarray*}
|(1-\frac{z}{z_i})(1-\frac{z}{\overline{z_i}})|&=&|1-z(\frac{1}{z_i}+\frac{1}{\overline{z_i}})+\frac{z^2}{|z_i|^2}|\\
&\leq&1+|z||\frac{1}{z_i}+\frac{1}{\overline{z_i}}|+\frac{|z|^2}{|z_i|^2}\\
&\leq&(1+|z||\frac{1}{z_i}+\frac{1}{\overline{z_i}}|)(1+\frac{|z|^2}{|z_i|^2}).
\end{eqnarray*}

Thus, if $z\in \CC$ we have
\begin{eqnarray*}
|\Phi (z)|\leq |e^{bz}|F_1(|z|)F_2(|z|)F_3(|z|)
\end{eqnarray*}
where
\begin{eqnarray*}
F_1(z)&=&|C|z^m\prod _{z_i\in\RR}(1+\frac{z}{|z_i|}),\\
F_2(z)&=&\prod _{z_i\in \CC ^+}(1+z|\frac{1}{z_i}+\frac{1}{\overline{z_i}}|),\\
F_3(z)&=&\prod _{z_i\in \CC ^+}(1+\frac{z^2}{|z_i|^2}).
\end{eqnarray*}

Since
\begin{eqnarray*}
\sum _{k}|Re \frac{1}{z_k}|<\infty ,
\end{eqnarray*}
we have that $F_1$ and $F_2$ are of minimal type. Moreover it is easy to see that
\begin{eqnarray*}
\lim _{R\rightarrow \infty}\frac{\log |F_1(R)|}{R}=\lim _{R\rightarrow \infty}\frac{\log |F_2(R)|}{R}=0.
\end{eqnarray*}

Now if we arranges $\{z_k\},~z_k\in \CC^+$ increasely by norms $0< |z_1|=|\overline{z_1}|\leq |z_2|=|\overline{z_2}|\leq |z_3|=|\overline{z_3}|...$ and put $\lambda _k=|z_k|$ then Theorem 1 of $\S 17.2$ in \cite{lev} gives
\begin{eqnarray*}
\lim _{n\rightarrow\infty}\frac{n}{\lambda _n}=\frac{d}{2\pi} = \frac{\sigma -\mu}{2\pi} .
\end{eqnarray*}
So we can apply Theorem 2 of $\S 12.1$ in \cite{lev} for $F_3(iz)$ to get
\begin{eqnarray*}
\lim _{R\rightarrow\infty}\frac{\log |F_3(R)|}{ R}=\frac{d}{2}= \frac{\sigma -\mu}{2}.
\end{eqnarray*}

Now if $R>0$ we have
\begin{eqnarray*}
\frac{\log |\Phi (R)|}{R}&\leq&b+\frac{\log |F_1(R)|}{R}+\frac{\log |F_2 (R)|}{R}+\frac{\log |F_3 (R)|}{R},\\
\frac{\log |\Phi (-R)|}{-R}&\geq&b+\frac{\log |F_1(R)|}{-R}+\frac{\log |F_2 (R)|}{-R}+\frac{\log |F_3 (R)|}{-R}.
\end{eqnarray*}

Hence
\begin{eqnarray*}
\sigma =\limsup _{R\rightarrow \infty}\frac{\log |\Phi (R)|}{R}&\leq&\lim _{R\rightarrow\infty }[b+\frac{\log |F_1(R)|}{R}+\frac{\log |F_2 (R)|}{R}+\frac{\log |F_3 (R)|}{R}]\\
&=&b+\frac{\sigma -\mu}{2},\\
\mu =\liminf _{R\rightarrow \infty}\frac{\log |\Phi (-R)|}{-R}&\geq&\lim _{R\rightarrow \infty}[b+\frac{\log |F_1(R)|}{-R}+\frac{\log |F_2 (R)|}{-R}+\frac{\log |F_3 (R)|}{-R}]\\
&=&b-\frac{\sigma -\mu}{2}.
\end{eqnarray*}

From above two inequalities we get
\begin{eqnarray*}
b=\frac{\sigma +\mu}{2}.
\end{eqnarray*}
\end{proof}

Now we consider the general case of kernel $\varphi \in L^1(\RR )$ satisfying condition (\ref{150506}).
\begin{lemma}
A function $\varphi\in L^1(\RR )$ satisfies condition (\ref{150506}) if and only if it satisfies condition (\ref{150506.1}). In this case $\Phi (z)$ defined by (\ref{fou}) is an entire function.
\label{lem4}\end{lemma}
\begin{proof}
If $\varphi$ satisfies condition (\ref{150506}), for each $n\in\NN$ there is $A_n<\infty$ such that
\begin{eqnarray*}
\int _{-\infty}^{\infty}e^{n|t|}|\varphi (t)|dt\leq A_n.
\end{eqnarray*}
Then
\begin{eqnarray*}
0&\geq& \limsup _{s\rightarrow\infty}\frac{\log \int _{|t|\geq s}e^{|t|n}|\varphi |dt}{s}\\
&\geq& \limsup _{s\rightarrow\infty}\frac{\log \int _{|t|\geq s}e^{sn}|\varphi |dt}{s}\\
&=&n+ \limsup _{s\rightarrow\infty}\frac{\log \int _{|t|\geq s}|\varphi |dt}{s}.
\end{eqnarray*}
Because $n$ is arbitrary, we have
\begin{eqnarray*}
\lim _{s\rightarrow\infty}\frac{\log \int _{|t|\geq s}|\varphi (t)|dt}{s}=-\infty .
\end{eqnarray*}
So $\varphi$ satisfies (\ref{150506.1}).

Conversely, if $\varphi$ satisfies (\ref{150506.1}), we take $p(s)$ be as in (\ref{170506.2}), that is
\begin{eqnarray*}
p(s)=-\log \int _{|t|\geq s}|\varphi (t)|dt ,~s\geq 0.
\end{eqnarray*}

Consider  the function
\begin{eqnarray*}
H(s)=\int _{-\infty}^{\infty}e^{|t|s}|\varphi (t)|dt.
\end{eqnarray*}

Because $\frac{d}{ds}e^{-p(s)}=-(|\varphi (s)|+|\varphi (-s)|)$ and $\varphi$ satisfies (\ref{150506.1}), using integration by parts we get
\begin{eqnarray*}
H(s)&=&-e^{ts-p(t)}|_{0}^{\infty}+s\int _{0}^{\infty}e^{ts-p(t)}dt\\
&=&||\varphi ||_{L^1(\RR )}+s\int _{0}^{\infty}e^{ts-p(t)}dt<\infty
\end{eqnarray*}
for any $s\geq 0$, so $\varphi$ satisfies (\ref{150506.1}).

When $\varphi$ satisfies (\ref{150506.1}), $\Phi (z)$ is clearly an entire function and
\begin{eqnarray*}
|\Phi (z)|\leq H(|z|).
\end{eqnarray*}
\end{proof}

As Lemma \ref{lem4} shows, the growth of $\Phi (z)$ is determinated by the function
\begin{equation}
G(s)=\int _{0}^{\infty}e^{ts-p(t)}dt.\label{170506.10}
\end{equation}

From now on we assume that $\varphi$ satisfies (\ref{150506}). For the sake of simplicity, we assume also that $\varphi$ has a non-compact support, so $p(s)<\infty$ for all $s\geq 0$.

It turns out that the growth of $G(s)$ in (\ref{170506.10}) is related to Young dual function $p^*(s)$ of $p(s)$, which is defined by
\begin{eqnarray*}
p^{*}(s)=\sup _{t\geq 0}[st-p(t)],~s\geq 0.
\end{eqnarray*}
It is easy to verify that $p^{*}$ is convex and satisfies
\begin{eqnarray*}
p(t)+p^*(s)&\geq&st,~s,t\geq 0,\\
\lim _{s\rightarrow\infty}\frac{p^{*}(s)}{s}&=&\infty .
\end{eqnarray*}
(see for e.g. Lecture 25 in \cite{lev}). So we can define Young dual function $p^{**}$ of $p^{*}$. This function satisfies $p^{**}(t)\leq p(t)$ for all $t\geq 0$.

We have the following result
\begin{theorem}
Let $\varphi$ be a function of $L^1(\RR )$ with non-compact support. Define p(s) as in (\ref{170506.2}) and Young dual function $p^{*}(s)$ and Young double-dual function $p^{**}(s)$ as above. Define $G(s)$ by (\ref{170506.10}). Then
\begin{eqnarray*}
\liminf _{s\rightarrow \infty}\frac{\log G(s)}{p^{*}(s)}\geq 1,
\end{eqnarray*}
and for any $\kappa >0$
\begin{eqnarray*}
\limsup _{s\rightarrow \infty}\frac{\log G(s)}{p^{*}(s+\kappa )}\leq 1.
\end{eqnarray*}
If additionally we have one of the following two conditions
\begin{equation}
\exp\{-p(t)+\frac{1}{\gamma}p(t\gamma )\}\in L^{1}(\RR ),\label{170506.11}
\end{equation}
for all $1>\gamma >0$ close enough to $1$, or
\begin{equation}
\exp\{-p^{**}(t)+\frac{1}{\gamma}p^{**}(t\gamma )\}\in L^{1}(\RR ),\label{170506.12}
\end{equation}
for all $1>\gamma >0$ close enough to $1$, then
\begin{eqnarray*}
\lim _{s\rightarrow\infty}\frac{\log G(s)}{p^{*}(s)}=1.
\end{eqnarray*}
\label{theo3}\end{theorem}
\begin{proof}
Because of properties of $p(s)$, for each $s>0$ large enough we can find a $t_s>0$ such that
\begin{eqnarray*}
st_s\geq p(t_s)+p^{*}(s)-\frac{1}{2},
\end{eqnarray*}
and
\begin{eqnarray*}
\lim _{s\rightarrow\infty}t_{s}=\infty .
\end{eqnarray*}
We have
\begin{eqnarray*}
G(s)=\int _{0}^{\infty}e^{ts-p(t)}dt\geq \int _{t_{s}-1}^{t_s}e^{ts-p(t)}dt.
\end{eqnarray*}
Now fixed $\gamma$ be any number less than $1$. By above inequality we have
\begin{eqnarray*}
e^{-\gamma p^{*}(s)}G(s)\geq e^{(1-\gamma)p^{*}(s)}\int _{t_{s}-1}^{t_s}e^{ts-p(t)-p^{*}(s)}dt.
\end{eqnarray*}
Because $p(t)$ in (\ref{170506.2}) is increasing, for $t_s-1\leq t\leq t_{s}$ we have
\begin{eqnarray*}
ts-p(t)-p^{*}(s)\geq s(t_s-1)-p(t_s)-p^{*}(s)\geq -s-\frac{1}{2}.
\end{eqnarray*}
So
\begin{eqnarray*}
e^{-\gamma p^{*}(s)}G(s)\geq e^{(1-\gamma)p^{*}(s)-s-\frac{1}{2}}.
\end{eqnarray*}
Because of
\begin{eqnarray*}
\lim _{s\rightarrow\infty}\frac{p^{*}(s)}{s}=\infty ,
\end{eqnarray*}
we get that
\begin{eqnarray*}
\liminf _{s\rightarrow\infty}\frac{\log G(s)}{p^{*}(s)}\geq 1.
\end{eqnarray*}
Now for any $\kappa >0$ we have $p(t)+p^{*}(s+\kappa )\geq st+\kappa t $. Hence
\begin{eqnarray*}
e^{-p^{*}(s+\kappa)}G(s)&=&\int _{0}^{\infty}e^{ts-p(t)-p^{*}(s+\kappa )}dt\\
&\leq&\int _{0}^{\infty}e^{-\kappa t}dt =\frac{1}{\kappa}.
\end{eqnarray*}
Thus
\begin{eqnarray*}
\limsup _{s\rightarrow\infty}\frac{\log G(s)}{p^{*}(s+\kappa )}\leq 1.
\end{eqnarray*}

Now assume that (\ref{170506.11}) or (\ref{170506.12}) is satisfied. Fixed $1>\gamma >0$ close enough to $1$. We have
\begin{eqnarray*}
\sup _{s\geq 0}[st-\frac{1}{\gamma}p^{*}(s)]=\frac{1}{\gamma}\sup _{s\geq 0}[s.(\gamma t)-p^{*}(s)]=\frac{1}{\gamma}p^{**}(\gamma t).
\end{eqnarray*}
Hence for all $t\geq 0$ we have
\begin{eqnarray*}
st-p(t)-\frac{1}{\gamma}p^{*}(s)\leq -p(t)+\frac{1}{\gamma}p^{**}(\gamma t).
\end{eqnarray*}

Because $p^{**}(t)\leq p(t)$ from above inequality we get
\begin{eqnarray*}
st-p(t)-\frac{1}{\gamma}p^{*}(s)\leq \min\{ -p(t)+\frac{1}{\gamma}p(\gamma t),-p^{**}(t)+\frac{1}{\gamma}p^{**}(\gamma t)\}.
\end{eqnarray*}
Hence
\begin{eqnarray*}
e^{st-p(t)-\frac{1}{\gamma}p^{*}(s)}\leq \min\{e^{ -p(t)+\frac{1}{\gamma}p(\gamma t)},e^{ -p^{**}(t)+\frac{1}{\gamma}p^{**}(\gamma t)}\},
\end{eqnarray*}
for all $s\geq 0$.

We have
\begin{eqnarray*}
\lim _{s\rightarrow\infty}e^{st-p(t)-\frac{1}{\gamma}p^{*}(s)}=0,
\end{eqnarray*}
for all $t\geq 0$, so we can apply Lebesgue's dominated convergence theorem to get
\begin{eqnarray*}
\limsup _{s\rightarrow\infty}\frac{\log G(s)}{p^{*}(s)}\leq \frac{1}{\gamma},
\end{eqnarray*}
for all $1>\gamma >0$ close enough to $1$, so it follows
\begin{eqnarray*}
\limsup _{s\rightarrow\infty}\frac{\log G(s)}{p^{*}(s)}\leq 1.
\end{eqnarray*}

Because we proved  before that
\begin{eqnarray*}
\liminf _{s\rightarrow\infty}\frac{\log G(s)}{p^{*}(s)}\geq 1,
\end{eqnarray*}
we get that
\begin{eqnarray*}
\lim _{s\rightarrow\infty}\frac{\log G(s)}{p^{*}(s)}= 1.
\end{eqnarray*}
\end{proof}
\begin{remark}
If $\varphi\geq 0$ then $\Phi (z)$ in (\ref{fou}) satisfies
\begin{eqnarray*}
\max _{|z|\leq s}|\Phi (z)|\geq \frac{1}{2}H(s),
\end{eqnarray*}
where $H(s)$ is the function defined in the proof of Lemma \ref{lem4}. So Theorem \ref{theo3} shows that $p^{*}(|z|)$ is an appropriate quantity to estimate the growth of $\log |\Phi (z)|$.
\end{remark}
\begin{remark}If $\varphi$ has the support in $[0,1]$ and
\begin{eqnarray*}
\sigma =\inf\{a\in [0,1]:~\varphi |_{[a,1]}=0~a.e.\},
\end{eqnarray*}
then it is easy to see that
\begin{eqnarray*}
\lim _{s\rightarrow \infty}\frac{p^{*}(s)}{\sigma s}=1,
\end{eqnarray*}
so by Theorem \ref{theo3} we obtain again the familiar result
\begin{eqnarray*}
\limsup _{|z|\rightarrow\infty}\frac{\log |\Phi (z)|}{|z|}\leq \sigma .
\end{eqnarray*}
\end{remark}
\begin{remark}
The class of functions $p(t)$ satisfying conditions (\ref{150506.1}) and (\ref{170506.11}) are fairly large. For example we can take $p (t)=t\log t$ or $p(t)=t^{\gamma}$ for any $\gamma >1$.
\end{remark}
\section{The size of small-valued sets}
In this section, we estimate the size of small-valued sets of Fourier transforms, i.e., estimate the Lebesgue measure of the sets
\begin{eqnarray*}
B_{\epsilon}=\{x\in \RR :~|\widehat{\varphi}(x)|\leq \epsilon ,~|x|\leq R_{\epsilon}\},
\end{eqnarray*}
where $\varphi$ is a nontrivial function of $L^1(\RR )$ satisfying condition (\ref{150506.1}), $R_{\epsilon}$ depends on small numbers $\epsilon$ and
\begin{eqnarray*}
\widehat{\varphi}(x )=\int _{-\infty}^{\infty}e^{-itx}\varphi (t)dt ,~x\in\RR .
\end{eqnarray*}

We will estimate these small-valued sets in terms of  the function $p(s)$ defined by (\ref{170506.2}) or in terms of its Young dual function $p^{*}(s)$.

First we obtain estimates in terms of $p(s)$. We define
$$\varphi ^{\epsilon}(t)=\left \{ \begin{array}{ll}\varphi (t)&\mbox{if }  |t|\leq s_{\epsilon} \\0&\mbox{if } |t|>s_{\epsilon}. \end{array}\right .$$
and
\begin{eqnarray*}
\Phi ^{\epsilon}(z)&=&\int _{-\infty}^{\infty}e^{zt}\varphi ^{\epsilon} (t)dt\equiv \int _{-s_{\epsilon}}^{s_{\epsilon}}e^{zt}\varphi (t)dt,\\
\widehat{\varphi ^{\epsilon}}(z)&=&\int _{-\infty}^{\infty}e^{-izt}\varphi ^{\epsilon} (t)dt\equiv \int _{-s_{\epsilon}}^{s_{\epsilon}}e^{-izt}\varphi (t)dt,
\end{eqnarray*}
where $s_{\epsilon}$ is defined by $(\ref{150506.2})$.

It is easy to see that $\Phi ^{\epsilon}(-iz)=\widehat{\varphi ^{\epsilon}}(z)$ are entire functions of order $1$  and type $\leq s_{\epsilon}$. More explicitly
\begin{eqnarray*}
|\widehat{\varphi ^{\epsilon}}(z)|\leq ||\varphi ||_{L^1(\RR )}e^{s_{\epsilon}|z|}.
\end{eqnarray*}

By the definition of $s_{\epsilon}$, we see that if $x\in \RR$ then
\begin{eqnarray*}
|\widehat{\varphi ^{\epsilon}}(x)-\widehat{\varphi }(x)|\leq \epsilon .
\end{eqnarray*}

First we estimate the size of small-valued sets of $\Phi ^{\epsilon}$. We have the following result

\begin{theorem} Let $s_{\epsilon},~\Phi ^{\epsilon}$,
be as above. Let $q>\frac{1}{2},~\beta >0$ be constants. For $\epsilon >0$ small enough we choose $R_{\epsilon}$ to satisfy
\begin{eqnarray*}
[(q+\frac{1}{2})\log R_{\epsilon}+\log (15e^3)][\log ||\varphi ||_{L^1(\RR )}+2es_{\epsilon}R_{\epsilon}]=-\log (\epsilon ^{\beta}+\epsilon).
\end{eqnarray*}

If $\epsilon$ small enough then
\begin{eqnarray*}
m(B_{\epsilon ^{\beta}+\epsilon})\leq R_{\epsilon}^{-q+\frac{1}{2}},
\end{eqnarray*}
where
\begin{eqnarray*}
B_{\epsilon ^{\beta}+\epsilon}=\{z\in\CC :~|\Phi ^{\epsilon}(z)|\leq \epsilon ^{\beta}+\epsilon ,~|z|\leq R_{\epsilon}\}.
\end{eqnarray*}
\label{theo1}\end{theorem}
\begin{proof} Because $||\varphi ||_{L^1(\RR)}\not= 0$, there exists $x_0\in \RR$ such that $\widehat{\varphi }(x_0)\not =0$. Then there exists a constant $C_0>0$ such that
\begin{eqnarray*}
|\widehat{\varphi ^{\epsilon}}(x_0)|\geq C_0,
\end{eqnarray*}
if $\epsilon$ is small enough. Changing variable if necessarily, we may assume that $|\Phi ^{\epsilon}(0)|\geq C_0$ if $\epsilon$ is small enough. If we choose
\begin{eqnarray*}
\eta _{\epsilon}=R_{\epsilon}^{-q-\frac{1}{2}},
\end{eqnarray*}
then by Theorem 4 of $\S 11.3$  in \cite{lev} if $|z|\leq R_{\epsilon}$ then
\begin{eqnarray*}
|\Phi ^{\epsilon}(z)|\geq \exp\{-[(q+\frac{1}{2})R_{\epsilon}+\log (15e^3)][\log ||\varphi ||_{L^1(\RR )}+s_{\epsilon}R_{\epsilon}]\}=\epsilon ^{\beta}+\epsilon ,
\end{eqnarray*}
except a set of disks whose sum of radii is less than
\begin{eqnarray*}
\eta _{\epsilon}R_{\epsilon}\equiv R_{\epsilon}^{-q+\frac{1}{2}}.
\end{eqnarray*}
\end{proof}

Because $|\widehat{\varphi^{\epsilon}} (x)|\leq \epsilon ^{\beta}+\epsilon$  if  $|\widehat{\varphi } (x)|\leq \epsilon^{\beta}$, by Theorem \ref{theo1} we have immediately
\begin{theorem} Let assumptions and choose $R_{\epsilon}$ as in Theorem \ref{theo1}. If $\epsilon$ small enough then
\begin{eqnarray*}
m(B_{\epsilon ^{\beta}})\leq R_{\epsilon}^{-q+\frac{1}{2}},
\end{eqnarray*}
where
\begin{eqnarray*}
B_{\epsilon ^{\beta}}=\{x\in\RR :~|\widehat{\varphi }(x)|\leq \epsilon ^{\beta},~|x|\leq R_{\epsilon}\}.
\end{eqnarray*}
\label{theo2}\end{theorem}

From Theorem \ref{theo3} we have
\begin{equation}
|\Phi (z)|\leq ||\varphi||_{L^1(\RR )}+|z|e^{p^{*}(|z|+1)},\label{170506.20}
\end{equation}
for all $z\in\CC$. So applying  Theorem 4 of $\S 11.3$  in \cite{lev} directly to $\Phi (z)$ we get the following result
\begin{theorem} Let $\varphi$ have a non-compact support.  Let $q>\frac{1}{2}$ be a constant. For $\epsilon >0$ small enough, define $R_{\epsilon}$ by
\begin{eqnarray*}
[(q+\frac{1}{2})R_{\epsilon}+\log (15e^3)][1+\log (R_{2\epsilon})+p^{*}(2R_{\epsilon}+1)]=-\log \epsilon .
\end{eqnarray*}
Then
\begin{eqnarray*}
m(B_{\epsilon})\leq R_{\epsilon}^{-q+\frac{1}{2}},
\end{eqnarray*}
where $m(B_{\epsilon})$ is the Lebesgue's measure of the set
\begin{eqnarray*}
B_{\epsilon}=\{z\in\CC :~|\Phi (z)|\leq \epsilon ,~|z|\leq R_{\epsilon}\}.
\end{eqnarray*}
Moreover we have
\begin{eqnarray*}
\lim _{\epsilon\rightarrow 0}\frac{\log p^{*}(2R_{\epsilon}+1)}{\log\log \frac{1}{\epsilon}}=1.
\end{eqnarray*}
\label{theo4}\end{theorem}
\begin{proof}
The conclusion about the size of $B_{\epsilon}$ is obtained similarly that of Theorems \ref{theo1} and \ref{theo2}. Because
\begin{eqnarray*}
\lim _{s\rightarrow\infty}\frac{p^{*}(s)}{s}=\infty ,
\end{eqnarray*}
in view of the definition of $R_{\epsilon}$ we get
\begin{eqnarray*}
\lim _{\epsilon\rightarrow 0}\frac{\log p^{*}(2R_{\epsilon}+1)}{\log\log \frac{1}{\epsilon}}=1.
\end{eqnarray*}
\end{proof}
\section{Proof of Theorem \ref{Tikhonov}}
Because $f_0,f_{\epsilon}\in L^2(\RR )$ we have
\begin{eqnarray*}
||f_0-f_{\epsilon}||^2_{L^2(\RR )}=\frac{1}{2\pi}||F(f_0)-F(f_{\epsilon})||^2_{L^2(\RR )},
\end{eqnarray*}
where $F(f_0)(\lambda )=\int _{-\infty}^{\infty}f_0(x)e^{-i\lambda x}dx=\widehat{f_0},~F(f_\epsilon )(\lambda )=\int _{-\infty}^{\infty}f_{\epsilon}(x)e^{-i\lambda x}dx=\widehat{f_{\epsilon}}$.

We have
\begin{eqnarray*}
||F(f_0)-F(f_{\epsilon})||^2_{L^2(\RR )}=\int _{-\infty}^{\infty}|\frac{\widehat{g_0}\overline{\widehat{\varphi _0}}}{|\widehat{\varphi _0}|^2}-\frac{\widehat{g_{\epsilon }}\overline{\widehat{\varphi _{\epsilon}}}}{\delta _{\epsilon}+|\widehat{\varphi _{\epsilon}}|^2}|^2.
\end{eqnarray*}

Now
\begin{eqnarray*}
|\frac{\widehat{g_0}\overline{\widehat{\varphi _0}}}{|\varphi _0|^2}-\frac{\widehat{g_{\epsilon }}\overline{\widehat{\varphi _{\epsilon}}}}{\delta _{\epsilon}+|\widehat{\varphi _{\epsilon}}|^2}|&\leq&|\frac{\widehat{g_0}\overline{\widehat{\varphi _0}}}{|\widehat{\varphi _0}|^2}-\frac{\widehat{g_{0}}\overline{\widehat{\varphi _{0}}}}{\delta _{\epsilon}+|\widehat{\varphi _{0}}|^2}|+|\frac{\widehat{g_0}\overline{\widehat{\varphi _0}}}{\delta _{\epsilon}+|\widehat{\varphi _0}|^2}-\frac{\widehat{g_{\epsilon }}\overline{\widehat{\varphi _{\epsilon}}}}{\delta _{\epsilon} +|\widehat{\varphi _{\epsilon}}|^2}|\\
&\leq&\frac{\delta _{\epsilon}|\widehat{g_0}||\overline{\widehat{\varphi _0}}|}{|\widehat{\varphi _0}|^2(\delta _{\epsilon}+|\widehat{\varphi _0}|^2)}+\frac{\delta _{\epsilon}|\widehat{g_0}\overline{\widehat{\varphi _0}}-\widehat{g_{\epsilon}}\overline{\widehat{\varphi _{\epsilon}}}|}{(\delta _{\epsilon}+|\widehat{\varphi _0}|^2)(\delta _{\epsilon}+|\widehat{\varphi _{\epsilon}}|^2)}\\
&&+\frac{|\widehat{\varphi _0}||\widehat{\varphi _{\epsilon}}||\widehat{g_0}\widehat{\varphi _0}-\widehat{g_{\epsilon}}\widehat{\varphi _{\epsilon}}|}{(\delta _{\epsilon}+|\widehat{\varphi _0}|^2)(\delta _{\epsilon}+|\widehat{\varphi _{\epsilon}}|^2)}.
\end{eqnarray*}

Now using Cauchy's inequality we have
\begin{eqnarray*}
\frac{\delta _{\epsilon}|\widehat{g_0}||\overline{\widehat{\varphi _0}}|}{|\widehat{\varphi _0}|^2(\delta _{\epsilon}+|\widehat{\varphi _0}|^2)}&\leq&\min \{\frac{\widehat{g_0}}{\widehat{\varphi _0}}=|F(f_0)|,\frac{\delta _{\epsilon}|\widehat{g_0}|}{|\widehat{\varphi _0}|^3}\},\\
\frac{\delta _{\epsilon}|\widehat{g_0}\overline{\widehat{\varphi _0}}-\widehat{g_{\epsilon}}\overline{\widehat{\varphi _{\epsilon}}}|}{(\delta _{\epsilon}+|\widehat{\varphi _0}|^2)(\delta _{\epsilon}+|\widehat{\varphi _{\epsilon}}|^2)}&\leq&\frac{|\widehat{g_0}\overline{\widehat{\varphi _0}}-\widehat{g_{\epsilon}}\overline{\widehat{\varphi _{\epsilon}}}|}{\delta _{\epsilon}}\\
\frac{|\widehat{\varphi _0}||\widehat{\varphi _{\epsilon}}||\widehat{g_0}\widehat{\varphi _0}-\widehat{g_{\epsilon}}\widehat{\varphi _{\epsilon}}|}{(\delta _{\epsilon}+|\widehat{\varphi _0}|^2)(\delta _{\epsilon}+|\widehat{\varphi _{\epsilon}}|^2)}&\leq&\frac{|\widehat{g_0}\widehat{\varphi _0}-\widehat{g_{\epsilon}}\widehat{\varphi _{\epsilon}}|}{\delta _{\epsilon}}.
\end{eqnarray*}

Hence
\begin{eqnarray*}
\int _{\RR }|F(f_0)-F(f_{\epsilon})|^2d\lambda &\leq&3[\int _{\RR}(\frac{\delta _{\epsilon}|\widehat{g_0}||\overline{\widehat{\varphi _0}}|}{|\widehat{\varphi _0}|^2(\delta _{\epsilon}+|\widehat{\varphi _0}|^2)})^2+\int _{\RR}(\frac{|\widehat{g_0}\widehat{\varphi _0}-\widehat{g_{\epsilon}}\widehat{\varphi _{\epsilon}}|}{\delta _{\epsilon}})^2\\
&&+\int _{\RR}(\frac{|\widehat{g_0}\overline{\widehat{\varphi _0}}-\widehat{g_{\epsilon}}\overline{\widehat{\varphi _{\epsilon}}}|}{\delta _{\epsilon}})^2]\\
&=&3[\int _{\RR}(\frac{\delta _{\epsilon}|\widehat{g_0}||\overline{\widehat{\varphi _0}}|}{|\widehat{\varphi _0}|^2(\delta _{\epsilon}+|\widehat{\varphi _0}|^2)})^2+\frac{1}{\delta _{\epsilon}^2}\int _{\RR}|\widehat{g_0}\widehat{\varphi _0}-\widehat{g_{\epsilon}}\widehat{\varphi _{\epsilon}}|^2\\
&&+\frac{1}{\delta _{\epsilon}^2}\int _{\RR}|\widehat{g_0}\overline{\widehat{\varphi _0}}-\widehat{g_{\epsilon}}\overline{\widehat{\varphi _{\epsilon}}}|^2].
\end{eqnarray*}

Now we write
\begin{eqnarray*}
&&\int _{\RR}(\frac{\delta _{\epsilon}|\widehat{g_0}||\overline{\widehat{\varphi _0}}|}{|\widehat{\varphi _0}|^2(\delta _{\epsilon}+|\widehat{\varphi _0}|^2)})^2\\
&=&\int _{|\lambda |>R _{\epsilon},|\widehat{\varphi _0} (\lambda )|<\epsilon ^{\beta}}(\frac{\delta _{\epsilon}|\widehat{g_0}||\overline{\widehat{\varphi _0}}|}{|\widehat{\varphi _0}|^2(\delta _{\epsilon}+|\widehat{\varphi _0}|^2)})^2+\int _{|\widehat{\varphi} (\lambda )|>\epsilon ^{\beta}}(\frac{\delta _{\epsilon}|\widehat{g_0}||\overline{\widehat{\varphi _0}}|}{|\widehat{\varphi _0}|^2(\delta _{\epsilon}+|\widehat{\varphi _0}|^2)})^2\\
&&+\int _{|\lambda |<R _{\epsilon},|\widehat{\varphi _0} (\lambda )|<\epsilon ^{\beta}}(\frac{\delta _{\epsilon}|\widehat{g_0}||\overline{\widehat{\varphi _0}}|}{|\widehat{\varphi _0}|^2(\delta _{\epsilon}+|\widehat{\varphi _0}|^2)})^2\\
&\leq&\int _{|\lambda |>R _{\epsilon},|\widehat{\varphi _0} (\lambda )|<\epsilon ^{\beta}}(\frac{\widehat{g_0}||\overline{\widehat{\varphi _0}}|}{|\widehat{\varphi _0}|^2})^2+\int _{|\widehat{\varphi _0} (i\lambda )|>\epsilon ^{\beta}}(\frac{\delta _{\epsilon}|\widehat{g_0}|}{|\widehat{\varphi _0}|^3})^2\\
&&+\int _{|\lambda |<R _{\epsilon},|\widehat{\varphi _0} (\lambda )|<\epsilon ^{\beta}}(\frac{\widehat{g_0}||\overline{\widehat{\varphi _0}}|}{|\widehat{\varphi _0}|^2})^2\\
&\leq& \int _{|\lambda |>R _{\epsilon},|\widehat{\varphi _0} (\lambda )|<\epsilon ^{\beta}}|F(f_0)|^2+\frac{\delta _{\epsilon}^2}{(\epsilon ^{3\beta})^2}\int _{|\widehat{\varphi _0} (\lambda )|>\epsilon ^{\beta}}|\widehat{g_0}|^2\\
&&+\int _{|\lambda |<R _{\epsilon},|\widehat{\varphi _0} (\lambda _{\epsilon})|<\epsilon ^{\beta}}|F(f_0)|^2\\
&\leq& \int _{|\lambda |>R _{\epsilon},|\widehat{\varphi _0} (\lambda )|<\epsilon ^{\beta}}|F(f_0)|^2+C_2\frac{\delta _{\epsilon}^2}{(\epsilon ^{3\beta})^2}+\int _{|\lambda |<R _{\epsilon},|\widehat{\varphi _0} (\lambda _{\epsilon})|<\epsilon ^{\beta}}|F(f_0)|^2.
\end{eqnarray*}

We have $||\widehat{\varphi _0}-\widehat{\varphi _{\epsilon}}||_{L^{\infty}(\RR)}\leq ||\varphi _0-\varphi _{\epsilon}||_{L^1(\RR )}$, a standard procedure give easily that $\lim _{\epsilon\rightarrow 0}||\widehat{g_0}\widehat{\varphi _0}-\widehat{g_{\epsilon}}\widehat{\varphi _{\epsilon}}||^2_{L^2(\RR )}=0$. Indeed, we have
\begin{eqnarray*}
||\widehat{g_0}\widehat{\varphi _0}-\widehat{g_{\epsilon}}\widehat{\varphi _{\epsilon}}||^2_{L^2(\RR )}&\leq& 2[||\widehat{\varphi _{\epsilon}}||^2_{L^{\infty}(\RR )}||\widehat{g_0}-\widehat{g_{\epsilon}}||^2_{L^2(\RR )}\\
&&+||\widehat{g_0}||^2_{L^2(\RR )}||\widehat{\varphi _0}-\varphi _{\epsilon}||^2_{L^{\infty}(\RR )}]\\
&\leq&2((||\varphi _0||_{L^1(\RR )}+\epsilon )^2+||g_0||_{L^2(\RR )}^2)\epsilon ^2\\
&\leq&2(||\varphi _0||_{L^1(\RR )}^2+||g_0||_{L^2(\RR )}^2+1)\epsilon ^2
\end{eqnarray*}

Similarly we have
\begin{eqnarray*}
||\widehat{g_0}\overline{\widehat{\varphi _0}}-\widehat{g_{\epsilon}}\overline{\widehat{\varphi _{\epsilon}}}||^2_{L^2(\RR )}\leq 2(||\varphi _0||_{L^1(\RR )}^2+||g_0||_{L^2(\RR )}^2+1)\epsilon ^2.
\end{eqnarray*}

So we have
\begin{eqnarray*}
\frac{1}{\delta _{\epsilon}^2}\int _{\RR}|\widehat{g_0}\widehat{\varphi _0}-\widehat{g_{\epsilon}}\widehat{\varphi _{\epsilon}}|^2+\frac{1}{\delta _{\epsilon}^2}\int _{\RR}|\widehat{g_0}\overline{\widehat{\varphi _0}}-\widehat{g_{\epsilon}}\overline{\widehat{\varphi _{\epsilon}}}|^2\leq C_1\frac{\epsilon ^2}{\delta _{\epsilon}^2}.
\end{eqnarray*}

Combining estimates above we get the conclusion of Theorem \ref{Tikhonov}.

\end{document}